\documentclass{amsart}

\usepackage{amsfonts}
\usepackage{amssymb}

\numberwithin{equation}{section}

\theoremstyle{plain}

\newtheorem{theorem}[subsection]{Theorem}
\newtheorem{proposition}[subsection]{Proposition}
\newtheorem{lemma}[subsection]{Lemma}
\newtheorem{claim}{Claim}

\theoremstyle{definition}

\newtheorem{problem}[subsection]{Problem}
\newtheorem*{example}{Example}

\renewcommand{\leq}{\leqslant}
\renewcommand{\geq}{\geqslant}

\providecommand{\supp}{\mathop{\rm supp}\nolimits}

\newcommand{\wt}{\widetilde}
\newcommand{\wh}{\widehat}

\newcommand{\Z}{\mathbb{Z}}

\newcommand{\R}{\mathbb{R}}
\newcommand{\C}{\mathbb{C}}

\begin{document}

\title{Three-term arithmetic progressions and sumsets}

\author{Tom Sanders}
\address{Department of Pure Mathematics and Mathematical Statistics\\
University of Cambridge\\
Wilberforce Road\\
Cambridge CB3 0WA\\
England } \email{t.sanders@dpmms.cam.ac.uk}

\begin{abstract}
Suppose that $G$ is an abelian group and $A \subset G$ is finite and
contains no non-trivial three-term arithmetic progressions. We show
that $|A+A| \gg_\varepsilon |A|(\log
|A|)^{\frac{1}{3}-\varepsilon}$.
\end{abstract}

\maketitle

\section{Introduction}

In \cite{GAF} Fre{\u\i}man proved the following qualitative
theorem.
\begin{theorem}[Fre{\u\i}man]
Suppose that $A \subset \Z$ is finite and contains no
non-trivial\footnote{A trivial three-term arithmetic progression is
one in which all three elements are the same.} three-term arithmetic
progressions. Then\footnote{By slight abuse of notation.} $|A+A|/|A|
\rightarrow \infty$ as $|A| \rightarrow \infty$.
\end{theorem}
The best known quantitative version of this theorem is achieved by
inserting Bourgain's most recent bound for Roth's theorem (see
\cite{JBRoth2}) into a result of Ruzsa's (see \cite{IZR3AP}).
\begin{theorem}[Bourgain-Ruzsa]\label{RuzBou}
Suppose that $A \subset \Z$ is finite and contains no non-trivial
three-term arithmetic progressions. Then
\begin{equation*}
|A+A| \gg |A| \left(\frac{\log |A|}{(\log \log
|A|)^3}\right)^{\frac{1}{6}}.
\end{equation*}
\end{theorem}
This theorem is interesting in its own right but has also been
applied (independently) by Schoen in \cite{TS} and Hegyv\'{a}ri,
Hennecart and Plagne in \cite{NHFHAP} to give a witty proof of the
following result regarding restricted sumsets.

If $A,B$ are subsets of an abelian group then we write
\begin{equation*}
A\ \wh{+}\ B:=\{a+b:a \in A, b \in B \textrm{ and } a \neq b\},
\end{equation*}
and call this the \emph{restricted sum} of $A$ and $B$.
\begin{theorem}[Schoen-Hegyv\'{a}ri-Hennecart-Plagne]\label{schoen}
Suppose $A$ and $B$ are two finite non-empty sets of integers, or
residues modulo an integer $m>1$, and put $n:=|A+B|$. Then
\begin{equation*}
\frac{|A\ \wh{+}\ B|}{|A+B|} = 1+O\left(\frac{(\log \log n)^3}{\log
n}\right)^{\frac{1}{6}}.
\end{equation*}
\end{theorem}
Recently a lot of work has been done on generalizing additive
problems in the integers to other abelian groups (see, for example,
\cite{BJGIZR,BJGTCTII,BJGTCTIII} and \cite{RM}) and in this paper we
not only improve the bounds in Theorems \ref{RuzBou} and
\ref{schoen} but we also extend them to cover arbitrary abelian
groups. Specifically our main result is the following theorem.
\begin{theorem}\label{premainthm}
Suppose that $G$ is an abelian group and $A \subset G$ is finite and
contains no non-trivial three-term arithmetic progressions. Then
\begin{equation*}
|A+A| \gg |A| \left(\frac{\log |A|}{(\log \log
|A|)^3}\right)^{\frac{1}{3}}.
\end{equation*}
\end{theorem}
This translates easily to an improvement of Theorem \ref{schoen}.
\begin{theorem}\label{schoent}
Suppose $A$ and $B$ are two finite non-empty subsets of an abelian
group $G$ and put $n:=|A+B|$. Then
\begin{equation*}
\frac{|A\ \wh{+}\ B|}{|A+B|} =1+O\left(\frac{(\log \log n)^3}{\log
n}\right)^{\frac{1}{3}}.
\end{equation*}
\end{theorem}
There are three main aspects to our arguments. First, to effect a
complete passage to general abelian groups we have to work slightly
harder when the sets in question have elements which differ by an
element of order 2. To deal with this we use a generalization of the
Bohr set technology of Bourgain \cite{JB}, as developed by Green and
the author in \cite{BJGTS2}.

Second, we use an energy increment argument in the style of
Heath-Brown \cite{DRHB} and Szemer{\'e}di \cite{ES} to prove a local
version of Roth's theorem which is particularly efficient
(essentially because of limitations in the modeling results of Green
and Ruzsa \cite{BJGIZR}) in our situation; this type of argument has
been deployed previously in \cite{TSASS}.

Finally we use a result which might be called a weak partially
polynomial version of the celebrated Fre{\u\i}man-Ruzsa theorem.
This type of result was first proved for finite fields by Green and
Tao in \cite{BJGTCTF}; the more general case we use was proved by
Green and the author in \cite{BJGTS2}.

The paper now splits into seven further sections. In
\S\S\ref{bohr}\verb!&!\ref{Ftools} we set up the basic machinery of
`local' Fourier analysis, which lets us prove our local version of
Roth's theorem in \S\ref{rbvar}. In \S\ref{bogchang} we prove the
partially polynomial version of the Fre{\u\i}man-Ruzsa theorem
before completing the main arguments in \S\ref{mainarg}.

In the final section, \S\ref{concluding remarks}, we discuss
improvements for particular groups $G$ and possible further
questions.

\section{Notation}

The book \cite{WR} serves as a general reference for the Fourier
transform, which we use throughout the paper.

Suppose that $G$ is a finite abelian group. $\wh{G}$ denotes the
\emph{dual group} of $G$, that is the group of homomorphisms
$\gamma:G\rightarrow S^1$, where $S^1:=\{z \in \C:|z|=1\}$, and we
write $M(G)$ for the space of measures on $G$ endowed with the norm
$\|.\|$ defined by $\|\mu\| := \int{d|\mu|}$.

There is one element of $M(G)$ worthy of particular note: the Haar
probability measure $\mu_G$. This measure is used to define the
Fourier transform which takes a function $f:G \rightarrow \C$ to
\begin{equation*}
\wh{f}:\wh{G} \rightarrow \C; \gamma \mapsto \int_{x \in
G}{f(x)\overline{\gamma}(x)d\mu_G(x)}=\frac{1}{|G|}\sum_{x \in
G}{f(x)\overline{\gamma}(x)}.
\end{equation*}
We use the Haar probability measure, $\mu_G$, on $G$ to define an
inner product on functions $f,g:G \rightarrow \C$ by
\begin{equation*}
\langle f,g\rangle:=\int_{x \in G}{f(x)\overline{g(x)}d\mu_G(x)}.
\end{equation*}
Since $\mu_G$ is normalized to be a probability measure,
Plancherel's theorem states that
\begin{equation*}
\langle f,g \rangle = \sum_{\gamma \in
\wh{G}}{\wh{f}(\gamma)\overline{\wh{g}(\gamma)}}.
\end{equation*}
Similarly we use $\mu_G$ to define the convolution of two
functions $f,g:G \rightarrow \C$:
\begin{equation*}
f \ast g(y):=\int_{x \in G}{f(y-x)g(x)d\mu_G(x)},
\end{equation*}
and a simple calculation tells us that $\wh{f \ast g}=\wh{f}\
\wh{g}$.

Finally it will sometimes be necessary to consider the Fourier
transform of a particularly complicated expression $E$. In this
case we may write $E^\wedge$ in place of $\wh{E}$.

\section{Bourgain systems}\label{bohr}

In \cite{JB}, Bourgain showed how to extend some of the techniques
of Fourier analysis from groups to a wider class of `approximate
groups'; in the paper \cite{BJGTS2} this was taken further when the
notion of a \emph{Bourgain system} was introduced. We refer the
reader to that paper for a more comprehensive discussion of Bourgain
systems and limit ourselves to recalling the key definitions and
tools that we shall require.

Suppose that $G$ is a finite abelian group and $d \geq 1$ is real. A
\emph{Bourgain system} $\mathcal{B}$ of dimension $d$ is a
collection $(B_\rho)_{\rho \in (0,2]}$ of subsets of $G$ such that
the following axioms are satisfied:
\begin{enumerate}
\item (Nesting) If $\rho' \leq \rho$ we have $B_{\rho'} \subseteq B_{\rho}$;
\item (Zero) $0 \in B_\rho$ for all $\rho \in (0,2]$;
\item (Symmetry) If $x \in B_{\rho}$ then $-x \in B_{\rho}$;
\item (Addition) For all $\rho,\rho'$ such that $\rho + \rho' \leq 1$ we have $B_{\rho} + B_{\rho'} \subseteq B_{\rho + \rho'}$;
\item (Doubling) If $\rho \leq 1$ then there is a set $X$ with $|X| \leq 2^d$
and \begin{equation*}B_{2\rho} \subset \bigcup_{x \in X}{x+B_\rho}.
\end{equation*}
\end{enumerate}
We define the $\emph{density}$ of $\mathcal{B}=(B_\rho)_\rho$ to be
$\mu_G(B_1)$ and denote it $\mu_G(\mathcal{B})$. Frequently we shall
consider several Bourgain systems
$\mathcal{B},\mathcal{B}',\mathcal{B}'',...$; in this case the
underlying sets will be denoted
$(B_\rho)_\rho,(B'_\rho)_\rho,(B''_\rho)_\rho,...$, and we shall
write $B,B',B'',...$ for the sets $B_1,B'_1,B''_1,...$.

\begin{example}[Bohr sets] There
is a natural valuation on $S^1$ defined by $\|z\|:=(2\pi)^{-1}|\arg
z|$, where $\arg$ is taken as mapping into $(-\pi,\pi]$. If $\Gamma
\subset \wh{G}$ and $\delta \in (0,1]$ then we put
\begin{equation*}
B(\Gamma,\delta):=\{x \in G: \|\gamma(x)\| \leq \delta \textrm{ for
all } \gamma \in \Gamma\},
\end{equation*}
and call such a set a \emph{Bohr set}.

It turns out that the system $(B(\Gamma,\rho\delta))_{\rho}$ is a
Bourgain system of density at least $\delta^{|\Gamma|}$ and
dimension $2|\Gamma|$, as the next lemma shows. By a slight abuse we
call this the Bourgain system \emph{induced} by the Bohr set
$B(\Gamma,\delta)$.
\begin{lemma}\label{bohrsize}
Suppose that $B(\Gamma,\delta)$ is a Bohr set. Then
\begin{equation*}
\mu_G(B(\Gamma,\delta)) \geq \delta^{|\Gamma|}
\end{equation*}
and there is a set $X$ of size at most $4^{|\Gamma|}$ such that
\begin{equation*}
B(\Gamma,2\delta) \subset \bigcup_{x \in X}{x+B(\Gamma,\delta)}.
\end{equation*}
\end{lemma}
The proof of this lemma is a simple averaging argument which may be
found, for example, in \cite[Lemma 4.20]{TCTVHV}.
\end{example}

Returning to Bourgain systems in general, we say that a Bourgain
system $\mathcal{B}'$ is a \emph{sub-system} of $\mathcal{B}''$ if
$B'_\rho \subset B''_\rho$ for all $\rho$. We shall be very
interested in sub-systems and consequently the following dilation
and intersection lemmas will be important. The first lemma is
immediate.
\begin{lemma}\label{bourgainsize}
Suppose that $\mathcal{B}$ is a Bourgain system of dimension $d$ and
$\lambda \in (0,1]$ is a parameter. Then $\lambda
\mathcal{B}:=(B_{\lambda\rho})_{\rho}$ is a Bourgain system of
dimension $d$ and density at least
$(\lambda/2)^d\mu_G(\mathcal{B})$.
\end{lemma}
\begin{lemma}\label{join}
Suppose that $\mathcal{B}^{(1)},\dots,\mathcal{B}^{(k)}$ are
Bourgain systems of dimensions $d_1,\dots,d_k$ respectively. Then
$\bigcap_{i=1}^k{\mathcal{B}^{(i)}}:=(\bigcap_{i=1}^{k}{B^{(i)}_\rho})_{\rho}$
is a Bourgain system of dimension at most $2(d_1+\dots+d_k)$ and
density at least
$4^{-(d_1+\dots+d_{k-1})}2^{-d_k}\prod_{i=1}^k{\mu_G(\mathcal{B}^{(i)})}$.
\end{lemma}
\begin{proof}
The conclusion is trivial apart from the doubling and density
estimates. For each $i$ with $1 \leq i \leq k$ there is a set $T_i$
with $|T_i| \leq 4^{d_i}$ such that $B^{(i)}_{2\rho} \subset T_i+
B^{(i)}_{\rho/2}$. Define a set $T$ as follows: for each
$(t_1,\dots,t_k) \in T_1 \times \dots \times T_k$ place one element
of $\bigcap_{i=1}^k{(t_i+B^{(i)}_{\rho/2})}$ in $T$ iff that set is
non-empty.

Now, if $t_0 \in \bigcap_{i=1}^k{(t_i+B^{(i)}_{\rho/2})}$ then the
map $t \mapsto t-t_0$ maps $\bigcap_{i=1}^k{(t_i+B^{(i)}_{\rho/2})}$
into $\bigcap_{i=1}^k{B^{(i)}_{\rho}}$, whence
\begin{equation*}
\bigcap_{i=1}^k{B^{(i)}_{2\rho}} \subset
T+\bigcap_{i=1}^k{B^{(i)}_{\rho}},
\end{equation*}
and the intersection has dimension at most $2(d_1+\dots+d_k)$.

The density estimate proceeds similarly. For each $i$ with $1 \leq i
\leq k-1$ let $T_i$ be a maximal subset of $G$ such that the sets
$(t+B_{1/4}^{(i)})_{t \in T_i}$ are disjoint. It follows that $|T|
\leq 4^{d_1}\mu_G(\mathcal{B}^{(i)})^{-1}$ and
\begin{equation*}
G \subset B_{1/4}^{(i)}-B_{1/4}^{(i)} + T_i \subset
B_{1/2}^{(i)}+T_i.
\end{equation*}
Thus there are some $x_1,\dots,x_{k-1} \in G$ such that
\begin{equation*}
\mu_G(\bigcap_{i=1}^{k-1}{(x_i+B_{1/2}^{(i)})}\cap B_{1/2}^{(k)})
\geq
4^{-(d_1+\dots+d_{k-1})}2^{-d_k}\prod_{i=1}^k{\mu_G(\mathcal{B}^{(i)})}.
\end{equation*}
Now for fixed $x_0 \in \bigcap_{i=1}^{k-1}{(x_i+B_{1/2}^{(i)})}\cap
B_{1/2}^{(k)}$ the map $x \mapsto x-x_0$ is an injection from
$\bigcap_{i=1}^{k-1}{(x_i+B_{1/2}^{(i)})}\cap B_{1/2}^{(k)}$ into
$\bigcap_{i=1}^k{B_1^{(i)}}$. The result follows.
\end{proof}
Not all Bourgain systems behave as regularly as we would like; we
say that a Bourgain system $\mathcal{B}$ of dimension $d$ is
\emph{regular} if
\begin{equation*}
1-2^3d|\eta| \leq \frac{\mu_G(B_1)}{\mu_G(B_{1+\eta})}\leq 1+
2^3d|\eta|
\end{equation*}
for all $\eta$ with $d|\eta| \leq 2^{-3}$. Typically, however,
Bourgain systems are regular, a fact implicit in the proof of the
following proposition.
\begin{proposition}\label{ubreg}
Suppose that $\mathcal{B}$ is a Bourgain system of dimension $d$.
Then there is a $\lambda \in [1/2,1)$ such that $\lambda\mathcal{B}$
is regular.
\end{proposition}
\begin{proof}
Let $f : [0,1] \rightarrow \R$ be the function $f(\alpha) :=
-\frac{1}{d}\log_2 \mu_G(B_{2^{-\alpha}})$ and note that $f$ is
non-decreasing in $\alpha$ with $f(1) - f(0) \leq 1$. We claim that
there is an $\alpha \in [\frac{1}{6}, \frac{5}{6}]$ such that
$|f(\alpha+ x) - f(\alpha)| \leq 3|x|$ for all $|x| \leq
\frac{1}{6}$. If no such $\alpha$ exists then for every $\alpha \in
[\frac{1}{6}, \frac{5}{6}]$ there is an interval $I(\alpha)$ of
length at most $\frac{1}{6}$ having one endpoint equal to $\alpha$
and with $\int_{I(\alpha)} df > \int_{I(\alpha)} 3 dx$. These
intervals cover $[\frac{1}{6}, \frac{5}{6}]$, which has total length
$\frac{2}{3}$. A simple covering lemma allows us to pass to a
disjoint subcollection $I_1 \cup ... \cup I_n$ of these intervals
with total length at least $\frac{1}{3}$. However we now have
\begin{equation*}
1 \geq \int^1_0 df \geq \sum_{i=1}^n \int_{I_i} df > \sum_{i = 1}^n
\int_{I_i} 3 \, dx \geq 1,
\end{equation*}
a contradiction. It follows that there is an $\alpha$ such that
$|f(\alpha+x) - f(\alpha)| \leq 3|x|$ for all $|x| \leq
\frac{1}{6}$. Setting $\lambda := 2^{-\alpha}$, it is easy to see
that
\begin{equation*}
(1+|\eta|)^{-3d} \leq \frac{\mu_G(B_{\lambda})}{\mu_G(B_{(1 +
\eta)\lambda})} \leq (1+|\eta|)^{3d}
\end{equation*}
whenever $|\eta| \leq 1/6$. But if $3d|\eta| \leq 1/2$ then
$(1+|\eta|)^{-3d} \leq 1+6d|\eta|$ and $(1+|\eta|)^{-3d} \geq
1-6d|\eta|$; it follows that $\lambda \mathcal{B}$ is a regular
Bourgain system.
\end{proof}

\section{Fourier analysis local to Bourgain systems}\label{Ftools}

Regular Bourgain systems are the `approximate groups' to which we
extend Fourier analysis; there is a natural candidate for
`approximate Haar measure' on $\mathcal{B}$: if $(B_\rho)_\rho$ is a
Bourgain system then we write $\beta_\rho$ for the normalized
counting measure on $B_\rho$ and simply $\beta$ for $\beta_1$. We
adopt similar conventions to before for the Bourgain systems
$\mathcal{B}',\mathcal{B}'',...$. It is worth noting that the
normalized measures introduced here are different from those in
\cite{BJGTS2} where positivity of the Fourier transform was also
desired.
\begin{lemma}[Approximate Haar measure]\label{haar}
Suppose that $\mathcal{B}$ is a regular Bourgain system of dimension
$d$. If $y \in B_{\eta}$ then $\|(y+\beta) - \beta\| \leq 2^4
d\eta$.
\end{lemma}
\begin{proof}
Note that $\supp{((y+\beta) - \beta)} \subset B_{1+\eta} \setminus
B_{1-\eta}$ whence
\begin{equation*}
\|(y+\beta) - \beta\| \leq \frac{\mu_G(B_{1+\eta} \setminus
B_{1-\eta})}{\mu_G(B_1)} \leq 2^4d\eta,
\end{equation*}
by regularity.
\end{proof}

The next two lemmas reflect two ways in which we commonly use the
property of regularity.

\begin{lemma}\label{cont}
Suppose that $\mathcal{B}$ is a regular Bourgain system of dimension
$d$. If $f:G \rightarrow \C$ then
\begin{equation*}
\|f \ast \beta - f \ast \beta(x)\|_{L^\infty(x+\beta_\eta)} \leq
2^4\|f\|_{L^\infty(\mu_G)}d\eta.
\end{equation*}
\end{lemma}
\begin{proof}
Note that
\begin{eqnarray*}
|f \ast \beta(x+y) - f \ast \beta(x)| & = & |f \ast
((-y+\beta)-\beta)(x)|\\ & \leq &
\|f\|_{L^\infty(\mu_G)}\|(-y+\beta)-\beta\|.
\end{eqnarray*}
The result follows by Lemma \ref{haar}.
\end{proof}

\begin{lemma}\label{nest}
Suppose that $\mathcal{B}$ is a regular Bourgain system of dimension
$d$ and $\kappa>0$ is a parameter. Then
\begin{equation*}
\{\gamma:|\wh{\beta}(\gamma)| \geq \kappa\} \subset
\{\gamma:|1-\gamma(x)| \leq 2^4d\kappa^{-1}\eta \textrm{ for all } x
\in B_\eta\}.
\end{equation*}
\end{lemma}
\begin{proof}
If $\gamma \in \{\gamma:|\wh{\beta}(\gamma)| \geq \kappa\}$ and $y
\in B_\eta$ then
\begin{eqnarray*}
\kappa |1-\gamma(y)| & \leq &
|\wh{\beta}(\gamma)||1-\overline{\gamma(y)}|\\ & = &
|\int{\gamma(x)d((y+\beta)-\beta)(x)}| \leq 2^4d\eta
\end{eqnarray*}
by Lemma \ref{haar}. The lemma follows.
\end{proof}

The final result of the section is a version of Bessel's inequality
local to Bourgain systems. Such a result was essentially proved by
Green and Tao in \cite[Corollary 8.6]{BJGTCTU3}, and serves to
replace some of the many applications of Parseval's theorem in the
local setting.

\begin{proposition}[Local Bessel inequality]\label{localbessel}
Suppose that $\mathcal{B}$ is a regular Bourgain system of dimension
$d$. Suppose that $f:G \rightarrow \C$ and $\epsilon \in (0,1]$ is a
parameter. Write $L_f:=\|f\|_{L^1(\beta)}^{-1}\|f\|_{L^2(\beta)}$.
Then there is a Bourgain system $\wt{\mathcal{B}}'$ of dimension
$2^2\epsilon^{-2}L_f^2$ such that
$\mathcal{B}':=\wt{\mathcal{B}}'\cap\mathcal{B}$ has
\begin{equation*}
\mu_G(\mathcal{B}') \geq
4^{-(d+2\epsilon^{-2}L_f^2)}\mu_G(\mathcal{B})
\end{equation*}
and
\begin{equation*}
\{\gamma :|\wh{fd\beta}(\gamma)| \geq \epsilon \|f\|_{L^1(\beta)}\}
\subset \{\gamma:|1-\gamma(x)| \leq 2^7(1+d)\epsilon^{-2}L_f^2\eta
\textrm{ for all } x \in B_{\eta}'\}.
\end{equation*}
\end{proposition}
To prove this we require an almost-orthogonality lemma due to
Cotlar \cite{MC}.
\begin{lemma}[Cotlar's almost orthogonality lemma]\label{cotlar} Suppose that
$v$ and $(w_j)$ are elements of an inner product space. Then
\begin{equation*}
\sum_j{|\langle v,w_j \rangle|^2} \leq \langle v ,v\rangle
\max_j{\sum_i{|\langle w_i,w_j\rangle|}}.
\end{equation*}
\end{lemma}
\begin{proof}
[Proof of Proposition \ref{localbessel}] Let
\begin{equation*}
S:=\{\gamma \in \wh{G}:|\wh{\beta}(\gamma)| \geq
\epsilon^2L_f^{-2}/2\},
\end{equation*}
and
\begin{equation*}
\Delta:=\{\gamma :|\wh{fd\beta}(\gamma)| \geq \epsilon
\|f\|_{L^1(\beta)}\}.
\end{equation*}
Pick  $\Lambda \subset \Delta$ maximal such that all the sets
$(\lambda+S)_{\lambda \in \Lambda}$ are disjoint. Now if $\gamma
\in \Delta$ then there is a $\lambda \in \Lambda$ such that
$\lambda+S \cap \gamma+S \neq \emptyset$ by maximality. It follows
that $\gamma \in \lambda +S-S$ i.e. $\Delta \subset \Lambda +
S-S$.

By Cotlar's lemma (Lemma \ref{cotlar}) we have
\begin{eqnarray*}
\sum_{\lambda \in \Lambda}{|\wh{fd\beta}(\lambda)|^2} &\leq &
\|f\|_{L^2(\beta)}^2\max_{\lambda \in \Lambda}{\sum_{\lambda' \in
\Lambda}{|\wh{\beta}(\lambda-\lambda')|}}\\ & \leq &
\|f\|_{L^2(\beta)}^2(1+|\Lambda|\epsilon^2L_f^{-2}/2),
\end{eqnarray*}
since $\lambda,\lambda' \in \Lambda$ and $\lambda-\lambda' \in S$
implies that $\lambda = \lambda'$. Since $\Lambda \subset \Delta$
we conclude that
\begin{equation*}
|\Lambda|\epsilon^2\|f\|_{L^1(\beta)}^2 \leq \sum_{\lambda \in
\Lambda}{|\wh{fd\beta}(\lambda)|^2}.
\end{equation*}
Combining all this we get that $ |\Lambda| \leq
2\epsilon^{-2}L_f^2$.

Let $\wt{\mathcal{B}}'$ be the Bourgain system induced by the Bohr
set $B(\Lambda,1)$ so $\mu_G(\wt{\mathcal{B}}') =1$ and $\dim
\wt{\mathcal{B}}' \leq 2|\Lambda| \leq 2^2\epsilon^{-2}L_f^2$.
Recalling that
\begin{equation*}
|1-\gamma(x)|=\sqrt{2(1-\cos (4\pi \|\gamma(x)\|))} \leq
4\pi\|\gamma(x)\|,
\end{equation*}
we certainly have
\begin{equation*}
\Lambda \subset \{\gamma:|1-\gamma(x)| \leq
2^6(1+d)\epsilon^{-2}L_f^2\eta \textrm{ for all } x \in
\wt{B}'_\eta\}.
\end{equation*}
By Lemma \ref{nest} $S$ is contained in
\begin{equation*}
\{\gamma:|1-\gamma(x)| \leq 2^5d\epsilon^{-2}L_f^2\eta \textrm{ for
all } x \in B_\eta\},
\end{equation*}
and so by the triangle inequality
\begin{equation*}
S-S \subset \{\gamma:|1-\gamma(x)| \leq 2^6d\epsilon^{-2}L_f^2\eta
\textrm{ for all } x \in B_\eta\}.
\end{equation*}
It follows that
\begin{equation*}
\Delta \subset \Lambda + S-S \subset \{\gamma:|1-\gamma(x)| \leq
2^7(1+d)\epsilon^{-2}L_f^2\eta \textrm{ for all } x \in B_\eta \cap
\wt{B}'_\eta\}.
\end{equation*}
The result follows by Lemma \ref{join} on letting
$\mathcal{B}':=\wt{\mathcal{B}}'\cap \mathcal{B}$.
\end{proof}

\section{A variant of the Bourgain-Roth theorem}\label{rbvar}

If $G$ is a finite group and $A \subset G$ then we can count the
number of three-term arithmetic progressions in $A$ using the
following trilinear form:
\begin{equation}\label{lamb}
\Lambda(f,g,h):=\int{f(x-y)g(x)h(x+y)d\mu_G(x)d\mu_G(y)}.
\end{equation}
This form has a well known Fourier expression gained by substituting
the inversion formul{\ae}  for $f,g$ and $h$ into (\ref{lamb}):
\begin{equation*}
\Lambda(f,g,h) = \sum_{\gamma \in \wh{G}
}{\wh{f}(\gamma)\wh{g}(-2\gamma)\wh{h}(\gamma)}.
\end{equation*}

In this section we shall prove the following result.
\begin{theorem}\label{bourgrothvar}
Suppose that $\mathcal{B}$ is a regular Bourgain system of dimension
$d$. Suppose that $A \subset G$ has $\alpha:=\|1_A \ast
\beta\|_{L^\infty(\mu_G)}$ -- that is the relative density of $A$ on
the translate of $B$ on which it is largest -- positive, and $A-A$
contains no elements of order 2. Then
\begin{equation*}
\Lambda(1_A,1_A,1_A) \geq
\left(\frac{\alpha}{2(1+d)}\right)^{2^{24}d\log \alpha^{-1} +
2^{52}\alpha^{-3}(\log \alpha^{-1})^2}\mu_G(\mathcal{B})^2.
\end{equation*}
\end{theorem}
We prove Theorem \ref{bourgrothvar} by iterating the following
lemma.
\begin{lemma}[Iteration lemma]\label{itlem}
Suppose that $\mathcal{B}$ is a regular Bourgain system of dimension
$d$. Suppose that $A \subset G$ has $\alpha:=\|1_A \ast
\beta\|_{L^\infty(\mu_G)}>0$ and $A-A$ contains no elements of order
2. Then at least one of the following is true.
\begin{enumerate}
\item \label{cse1} \emph{(Lots of three-term progressions.)}
\begin{equation*}
\Lambda(1_A,1_A,1_A) \geq
\frac{\alpha^3}{2^5}\left(\frac{\alpha^3}{2^{44}(1+d)^3}\right)^d\mu_G(\mathcal{B})^2.
\end{equation*}
\item \label{cse2} \emph{(Density increment I)} There is a regular
dilate $\mathcal{B}''$ of $\mathcal{B}$ with
\begin{equation*}
\mu_G(\mathcal{B}'') \geq
\left(\frac{\alpha^2}{2^{25}(1+d)^2}\right)^d\mu_G(\mathcal{B})
\end{equation*}
such that $\|1_A \ast \beta''\|_{L^\infty(\mu_G)} \geq
\alpha(1+2^{-12})$.
\item \label{wrd} \emph{(Density increment II)} There is a
regular dilate $\mathcal{B}'''$ of $(\{2x:x \in B_\rho\})_\rho$ with
\begin{equation*}
\mu_G(\mathcal{B}''') \geq
\frac{\alpha}{2^2}\left(\frac{\alpha^3}{2^{36}(1+d)^3}\right)^d\mu_G(\mathcal{B})
\end{equation*}
such that $\|1_A \ast \beta'''\|_{L^\infty(\mu_G)} \geq
\alpha(1+2^{-8})$.
\item \label{cse3} \emph{(Density increment III)} There is a
Bourgain system $\wt{\mathcal{B}}''''$ of dimension at most
$2^{13}\alpha^{-3}$ and a dilate $\mathcal{B}'''$ of $(\{2x:x \in
B_\rho\})_\rho$ such that their intersection, $\mathcal{B}''''$, is
regular with
\begin{equation*}
\mu_G(\mathcal{B}'''') \geq
\frac{\alpha}{2^2}\left(\frac{\alpha^3}{2^{22}(1+d)}\right)^{2^{13}\alpha^{-3}}
\left(\frac{\alpha^5}{2^{48}(1+d)^3}\right)^d\mu_G(\mathcal{B})
\end{equation*}
such that $\|1_A \ast \beta''''\|_{L^\infty(\mu_G)} \geq
\alpha(1+2^{-8})$.
\end{enumerate}
\end{lemma}
The different cases (\ref{cse2}), (\ref{wrd}) and (\ref{cse3}) are
the outcomes of different parts of the proof; we separate them for
ease of understanding.

The proof of the lemma requires the following technical result which
converts energy on non-trivial Fourier modes into a density
increment.
\begin{lemma}[$\ell^2$-density increment lemma]\label{l2increment}
Suppose that $\mathcal{B}$ is a regular Bourgain system of dimension
$d$. Suppose that $A \subset G$ has $\alpha:=1_A \ast \beta(0_G)>0$
and $c >0$ is a parameter. Write $\eta:=c\alpha/2^{10}(1+d)$ and
suppose that $\mathcal{B}'$ is a sub-system of $\eta\mathcal{B}$ and
that there is a set of characters
\begin{equation*}
\Lambda:=\{\gamma :|1-\gamma(x)| \leq 1/2 \textrm{ for all } x \in
B'\}
\end{equation*}
such that
\begin{equation*}
\sum_{\lambda \in \Lambda}{|((1_A-\alpha)1_{B})^\wedge(\lambda)|^2}
\geq c\alpha^2\mu_G(B).
\end{equation*}
Then $\|1_A \ast \beta'\|_{L^\infty(\mu_G)} \geq\alpha(1+c/2^3)$.
\end{lemma}
\begin{proof}
Write $f:=1_A - \alpha$. The triangle inequality shows that if
$\lambda \in \Lambda$ then
\begin{equation*}
|\wh{\beta'}(\lambda)| \geq \int{d\beta'} - \int{|1-\lambda|d\beta'}
\geq 1/2,
\end{equation*}
whereupon (from the hypothesis of the lemma)
\begin{equation*}
c\alpha^2\mu_G(B)/2^2 \leq \sum_{\gamma \in
\wh{G}}{|\wh{f1_{B}}(\gamma)\wh{\beta'}(\gamma)|^2}.
\end{equation*}
Plancherel's theorem (and dividing by $\mu_G(B)$) then gives
\begin{equation*}
\langle (f1_{B}) \ast \beta',(fd\beta)\ast \beta' \rangle \geq c
\alpha^2/2^2.
\end{equation*}
We expand this inner product:
\begin{eqnarray}
\label{dec} \langle (f1_{B}) \ast \beta',(fd\beta) \ast \beta'
\rangle & = & \langle (1_A1_{B}) \ast \beta', (1_Ad\beta) \ast
\beta' \rangle\\ \nonumber & & - \alpha\langle 1_{B} \ast \beta',
(1_Ad\beta)\ast \beta'\rangle
\\ \nonumber  & & - \alpha \langle (1_A1_{B}) \ast \beta', \beta \ast \beta'
\rangle\\ \nonumber & & + \alpha^2\langle 1_{B} \ast \beta', \beta
\ast \beta' \rangle.
\end{eqnarray}
We estimate the last three-terms: By Lemma \ref{haar} we have
\begin{eqnarray}
\label{er} \|\beta \ast \beta' \ast \beta' - \beta\| & \leq &
\int{\|(y+\beta) - \beta\|d(\beta' \ast \beta')(y)}\\ \nonumber
 & \leq & \sup_{y \in \supp \beta' \ast \beta'}{\|(y+\beta)
- \beta\|}\\ \nonumber & \leq & \sup_{y \in B_2'}{\|(y+\beta) -
\beta\|}\\ \nonumber & \leq & \sup_{y \in B_{2\eta}}{\|(y+\beta) -
\beta\|} \leq c\alpha/2^5.
\end{eqnarray}
Now
\begin{equation*}
\langle 1_{B} \ast \beta', (1_Ad\beta)\ast \beta'\rangle = \langle
\beta \ast \beta' \ast \beta', (1_A1_{B}) \rangle
\end{equation*}
and
\begin{equation*}
|\langle \beta \ast \beta' \ast \beta', 1_A1_{B} \rangle - \langle
\beta,1_A1_{B}\rangle| \leq c\alpha/2^5
\end{equation*}
by (\ref{er}); $\langle \beta,1_A1_{B}\rangle = \alpha$, so
\begin{equation*}
|\langle 1_{B} \ast \beta', (1_Ad\beta)\ast \beta'\rangle - \alpha|
\leq c\alpha/2^5.
\end{equation*}
By symmetry
\begin{equation*}
|\langle (1_A1_{B}) \ast \beta', \beta \ast \beta'\rangle - \alpha|
\leq c\alpha/2^5,
\end{equation*}
and similarly
\begin{equation*}
|\langle 1_{B} \ast \beta', \beta \ast \beta'\rangle - 1| \leq
c\alpha/2^5.
\end{equation*}
Inserting these last three estimates into (\ref{dec}) we get
\begin{equation*}
\langle (f1_{B}) \ast \beta',(fd\beta) \ast \beta' \rangle \leq
\langle (1_A1_{B}) \ast \beta', (1_Ad\beta) \ast \beta' \rangle -
\alpha^2 + c\alpha^2/2^3.
\end{equation*}
We conclude that
\begin{equation*}
\alpha^2(1+c/2^3) \leq \langle (1_A1_{B}) \ast \beta', (1_Ad\beta)
\ast \beta' \rangle.
\end{equation*}
Finally
\begin{eqnarray*}
\langle (1_A1_{B}) \ast \beta', (1_Ad\beta) \ast \beta' \rangle &
\leq &
\|(1_A1_{B}) \ast \beta'\|_{L^\infty(\mu_G)} \|(1_Ad\beta) \ast \beta'\|\\
& \leq &\|(1_A1_{B}) \ast
\beta'\|_{L^\infty(\mu_G)}\|1_A\|_{L^1(\beta)}\|\beta'\|\\ & = &
\|(1_A1_{B}) \ast \beta'\|_{L^\infty(\mu_G)}\alpha\\ & \leq & \|1_A
\ast \beta'\|_{L^\infty(\mu_G)}\alpha;
\end{eqnarray*}
we get the result on dividing by $\alpha$.
\end{proof}

\begin{proof}[Proof of Lemma \ref{itlem}]
Suppose that we are not in case (\ref{cse2}) of the lemma, so we may
certainly assume that for all regular dilates $\mathcal{B}''$ of
$\mathcal{B}$ with
\begin{equation*}
\mu_G(\mathcal{B}'') \geq
\left(\frac{\alpha^2}{2^{25}(1+d)^2}\right)^d\mu_G(\mathcal{B}),
\end{equation*}
we have
\begin{equation}\label{asmptn}
\|1_A \ast \beta''\|_{L^\infty(\mu_G)} \leq \alpha(1+2^{-12}).
\end{equation}

Apply Proposition \ref{ubreg} to pick $\lambda'$ so that
$\mathcal{B}':=\lambda'\mathcal{B}$ is regular and
\begin{equation*}
\alpha/2^{16}(1+d) \leq \lambda' < \alpha/2^{15}(1+d).
\end{equation*}
Apply Proposition \ref{ubreg} to pick $\lambda''$ so that
$\mathcal{B}'':=\lambda''\mathcal{B}'$ is regular and
\begin{equation*}
\alpha/2^8(1+d) \leq \lambda'' < \alpha/2^7(1+d).
\end{equation*}

Suppose that $\lambda \in [\lambda''\lambda',\lambda']$. A trivial
instance of Young's inequality tells us that
\begin{eqnarray*}\|1_A \ast \beta \ast \beta_{\lambda} -
1_A \ast \beta\|_{L^\infty(\mu_G)} & \leq &
\|1_A\|_{L^\infty(\mu_G)}\|\beta \ast \beta_{\lambda} - \beta\|\\ &
\leq &\int{\|(y+\beta) - \beta\|d\beta_{\lambda}(y)}\\ & \leq &
\sup_{y \in B_\lambda}{\|(y+\beta) - \beta\|}\\& \leq & 2^4d\lambda
\leq \alpha/2^{11}
\end{eqnarray*}
by Lemma \ref{haar} and the fact that $\lambda \leq \lambda'$. Let
$x' \in G$ be such that $1_A \ast \beta(x')=\alpha$. It follows from
the previous calculation that
\begin{equation*}
|(1_A \ast \beta_{\lambda} -\alpha)\ast \beta(x')| \leq
\alpha/2^{11}.
\end{equation*}
Moreover by assumption (\ref{asmptn}) (applicable by Lemma
\ref{bourgainsize} and the fact that $\lambda \geq
\lambda''\lambda'$) we have
\begin{equation*}
1_A \ast \beta_{\lambda} - \alpha \leq \alpha/2^{12}.
\end{equation*}
For functions $g:G \rightarrow \C$ we write $g_+:=(|g|+g)/2$ and
$g_-:=(|g|-g)/2=g_+-g$. Now, combining our last two expressions then
yields
\begin{eqnarray*}
|1_A \ast \beta_{\lambda}-\alpha| \ast \beta(x') &= & (1_A \ast
\beta_{\lambda}-\alpha)_+ \ast \beta(x')\\
& & + (1_A \ast \beta_{\lambda}-\alpha)_- \ast \beta(x')\\
& = & 2 (1_A \ast
\beta_{\lambda}-\alpha)_+ \ast \beta(x')\\
& & - (1_A \ast \beta_{\lambda}-\alpha) \ast \beta(x')\\
& \leq & \alpha/2^{10}.
\end{eqnarray*}
Applying this expression with $\lambda$ equals $\lambda'$ and
$\lambda''\lambda'$ we get
\begin{eqnarray*}
\alpha/2^9 & \geq & \left(|1_A \ast \beta'-\alpha| + |1_A \ast
\beta'' - \alpha|\right)\ast \beta(x')\\ & \geq & \inf_{x \in
G}{\left(|1_A \ast \beta'(x)-\alpha|+|1_A \ast \beta''(x)-\alpha|
\right)}.
\end{eqnarray*}
By translating $A$ we may assume that the infimum on the right is
attained at $x=0_G$; we write
\begin{equation*}
\alpha':=1_A \ast \beta'(0_G), \alpha'':=1_A \ast \beta''(0_G),
f':=1_A - \alpha', \textrm{ and } f'':=1_A - \alpha'',
\end{equation*}
and note that
\begin{equation*}
|\alpha''-\alpha| \leq \alpha/2^9 \textrm{ and } |\alpha'-\alpha|
\leq \alpha/2^9.
\end{equation*}
Now by trilinearity of $\Lambda$ we have
\begin{eqnarray}
\label{spt} \Lambda(1_A1_{B'},1_A1_{B''},1_A 1_{B'})& = &
\Lambda(1_A1_{B'},1_A1_{B''},\alpha'1_{B'})\\ \nonumber & & +
\Lambda(\alpha'1_{B''},1_A1_{B''},f'1_{B'})\\ \nonumber & & +
\Lambda(f'1_{B'},\alpha''1_{B''},f'1_{B'})\\ \nonumber & &
+\Lambda(f'1_{B'},f''1_{B''},f'1_{B'}).
\end{eqnarray}
We can easily estimate the first two terms on the right using the
following fact.
\begin{claim}
Suppose that $g:G \rightarrow \C$ has $\|g\|_{L^\infty(\mu_G)}
\leq 1$. Then
\begin{equation*}
|\Lambda(g1_{B'},1_A1_{B''},1_{B'}) - \alpha''g \ast
\beta'(0_G)\mu_G(B'')\mu_G(B')| \leq
\alpha''\alpha'\mu_G(B'')\mu_G(B')/2^2.
\end{equation*}
\end{claim}
\begin{proof}
Recall that $\Lambda(g1_{B'},1_A1_{B''},1_{B'})$ equals
\begin{equation*}
\int{g(x-y)1_{B'}(x-y)1_A(x)1_{B''}(x)1_{B'}(x+y)d\mu_G(x)d\mu_G(y)}
\end{equation*}
by definition. By the change of variables $u=x-y$ and symmetry of
$B'$ we conclude that this expression is in turn equal to
\begin{equation*}
\int{g(u)1_{B'}(u)1_A(x)1_{B''}(x)1_{B'}(u-2x)d\mu_G(x)d\mu_G(u)}.
\end{equation*}
Now the difference between this term and
\begin{equation*}
\int{g(u)1_{B'}(u)1_A(x)1_{B''}(x)1_{B'}(u)d\mu_G(u)d\mu_G(x)}(=\alpha''g
\ast \beta'(0_G)\mu_G(B'')\mu_G(B'))
\end{equation*}
is at most
\begin{equation*}
\|g\|_{L^\infty(\mu_G)}\int{1_{B'}(u)1_A(x)1_{B''}(x)|1_{B'}(u-2x) -
1_{B'}(u)|d\mu_G(x)d\mu_G(u)}
\end{equation*}
in absolute value. But if $x \in B''$ then $2x \in B''_1+B''_1
\subset B''_2=B'_{2\lambda''}$ whence if $u \in B'_{1-2\lambda''}$
then $1_{B'}(u) = 1_{B'}(u-2x)$. It follows that this error term is
at most
\begin{equation*}
\alpha''\mu_G(B'')\mu_G(B'_1 \setminus B'_{1-2\lambda''}) \leq
2^4\alpha''d\lambda''\mu_G(B'')\mu_G(B')
\end{equation*}
by regularity of $\mathcal{B}'$. The claim follows in view of the
earlier choice of $\lambda''$ and the fact that $\alpha' \geq
\alpha/2$.
\end{proof}
It follows by applying this claim with $g=1_A$ that
\begin{equation}\label{kb1}
|\Lambda(1_A1_{B'},1_A1_{B''},\alpha'1_{B'}) -
\alpha''\alpha'^2\mu_G(B'')\mu_G(B')| \leq
\alpha''\alpha'^2\mu_G(B'')\mu_G(B')/2^2.
\end{equation}
Moreover, since $f' \ast \beta'(0_G) = 0$ the claim applied with
$g=f'$ gives
\begin{equation}\label{kb2}
|\Lambda(\alpha'1_{B'},1_A 1_{B''},f'1_{B'})| \leq
\alpha''\alpha'^2\mu_G(B'')\mu_G(B')/2^2.
\end{equation}
In view of (\ref{kb1}), (\ref{kb2}) and the decomposition
(\ref{spt}) we conclude (by the triangle inequality) that either
\begin{enumerate}
\item \begin{equation*} |\Lambda(1_A1_{B'},1_A1_{B''},1_A 1_{B'})| \geq
\alpha''\alpha'^2\mu_G(B'')\mu_G(B')/2^2,
\end{equation*}
and we are in case (\ref{cse1}) of the lemma; \item or
\begin{equation*}
|\Lambda(f'1_{B'},\alpha''1_{B''},f'1_{B'})| \geq
\alpha''\alpha'^2\mu_G(B'')\mu_G(B')/2^3,
\end{equation*}
and it turns out that we are in case (\ref{wrd}) of the lemma;
\item or
\begin{equation*} |\Lambda(f'1_{B'},f''1_{B''},f'1_{B'})|\geq
\alpha''\alpha'^2\mu_G(B'')\mu_G(B')/2^3,
\end{equation*}
and it turns out that we are in case (\ref{cse3}) of the lemma.
\end{enumerate}
The first conclusion is immediate. The second and third are verified
(respectively) in the following two claims.
\begin{claim}
If
\begin{equation*}
|\Lambda(f'1_{B'},\alpha''1_{B''},f'1_{B'})| \geq
\alpha''\alpha'^2\mu_G(B'')\mu_G(B')/2^3
\end{equation*}
then we are in case (\ref{wrd}) of the lemma.
\end{claim}
\begin{proof}
In view of the Fourier expression for $\Lambda$ we get
\begin{equation}\label{coret}
\alpha'^2\mu_G(B'')\mu_G(B')/2^3  \leq \sum_{\gamma \in
\wh{G}}{|\wh{1_{B''}}(2\gamma)||\wh{f'1_{B'}}(\gamma)|^2}.
\end{equation}
It turns out that the characters for which $|\wh{1_{B''}}(2\gamma)|$
is large support a lot of the mass of the sum on the right: Let
$\epsilon =\alpha'/2^4$ and put
\begin{equation*}
\Lambda:=\{\gamma \in \wh{G}:|\wh{1_{B''}}(2\gamma)| \geq
\epsilon\mu_G(B'')\}.
\end{equation*}
Then
\begin{eqnarray*}
\sum_{\gamma \not\in
\Lambda}{|\wh{1_{B''}}(2\gamma)||\wh{f'1_{B'}}(\gamma)|^2} & \leq &
 \epsilon\mu_G(B'')\sum_{\gamma \in \wh{G}}{|\wh{f'1_{B'}}(\gamma)|^2}\\ & = & \epsilon\mu_G(B'')\mu_G(B')\|f'\|_{L^2(\beta')}^2,
\end{eqnarray*}
by the triangle inequality and Parseval's theorem. Now $
\|f'\|_{L^2(\beta')}^2 = \alpha' - \alpha'^2$, so it follows that
this last expression is at most $\epsilon\alpha'\mu_G(B'')\mu_G(B')$
and hence by the triangle inequality and (\ref{coret}) we have
\begin{equation*}
\alpha'^2\mu_G(B')/2^4 \leq \sum_{\gamma \in
\Lambda}{|\wh{f'1_{B'}}(\gamma)|^2}.
\end{equation*}
Note that $(\{2x:x \in B''_{\rho}\})_{\rho}$ is a Bourgain system of
dimension $d$. Apply Proposition \ref{ubreg} to pick $\lambda'''$ so
that $\mathcal{B}''':=\lambda'''(\{2x:x \in B''_{\rho}\})_\rho$ is
regular and
\begin{equation*}
\alpha/2^{11}(1+d) \leq \lambda''' < \alpha/2^{10}(1+d);
\end{equation*}
since $\mathcal{B}''$ is a dilate of $\mathcal{B}$, $\mathcal{B}'''$
is a dilate of $(\{2x:x \in B_\rho\})_{\rho}$. By Lemma \ref{nest}
we have that
\begin{eqnarray*}
\Lambda & \subset & \{\gamma : |1-(2\gamma)(x)| \leq 1/2 \textrm{
for all } x \in B''_{\lambda'''}\}\\ & = & \{\gamma : |1-\gamma(x)|
\leq 1/2 \textrm{ for all } x \in B'''\}.
\end{eqnarray*}
Now $\mathcal{B}'''$ is a subsystem of
$(\alpha'/2^{14}(1+d))\mathcal{B}'$ so we apply Lemma
\ref{l2increment} with $c=2^{-4}$ to see that
\begin{equation*}
\|1_A \ast \beta'''\|_{L^\infty(\mu_G)} \geq \alpha'(1+2^{-7}) \geq
\alpha(1+2^{-8}).
\end{equation*}
It remains only to verify the bound on the density of
$\mathcal{B}'''$. Note that
\begin{equation*}
\|1_A \ast \beta''_{\lambda'''} \ast \beta' - 1_A \ast
\beta'\|_{L^\infty(\mu_G)} \leq 2^4d\lambda'''\lambda'' \leq
\alpha'/2
\end{equation*}
by Lemma \ref{haar}. Whence
\begin{equation*}
1_A \ast \beta''_{\lambda'''}\ast \beta'(x') \geq 1_A \ast
\beta'(x') - \alpha'/2 \geq \alpha/2^2.
\end{equation*}
By averaging it follows that there is some $x'' \in G$ such that $
1_A \ast \beta''_{\lambda'''}(x'') \geq \alpha/2^2$. Since $A - A$
contains no elements of order 2 we have that $x \mapsto 2x$ is
injective when restricted to $A$; we conclude that
\begin{eqnarray*}
\mu_G(2B''_{\lambda'''}) & =& \mu_G(2(x''+B''_{\lambda'''}))\\
& \geq & \mu_G(2(A \cap(x''+B''_{\lambda'''})))\\& = &\mu_G(A
\cap(x''+B''_{\lambda'''}))\\ & \geq & \frac{\alpha}{2^2}
\mu_G(B''_{\lambda'''})\\ & \geq &
\frac{\alpha}{2^2}\left(\frac{\lambda'\lambda''\lambda'''}{2}\right)^d\mu_G(\mathcal{B}),
\end{eqnarray*}
by Lemma \ref{bourgainsize}. The claim follows.
\end{proof}

\begin{claim}
If
\begin{equation*}
|\Lambda(f'1_{B'},f''1_{B''},f'1_{B'})| \geq
\alpha''\alpha'^2\mu_G(B'')\mu_G(B')/2^3
\end{equation*}
then we are in case (\ref{cse3}) of the lemma.
\end{claim}
\begin{proof}
In view of the Fourier expression for $\Lambda$ we have
\begin{equation}\label{coret2}
\alpha''\alpha'^2\mu_G(B'')\mu_G(B')/2^3 \leq \sum_{\gamma \in
\wh{G}}{|\wh{f''1_{B''}}(2\gamma)||\wh{f1_{B'}}(\gamma)|^2}.
\end{equation}
As in the previous claim we may ignore the characters supporting
small values of $\wh{f''1_{B''}}(\gamma)$: Let $\epsilon
=\alpha''\alpha'/2^4$ and put
\begin{equation*}
\Lambda:=\{\gamma \in \wh{G}:|\wh{f''1_{B''}}(2\gamma)| \geq
\epsilon\mu_G(B'')\}.
\end{equation*}
Then
\begin{eqnarray*}
\sum_{\gamma \not\in
\Lambda}{|\wh{f''1_{B''}}(2\gamma)||\wh{f'1_{B'}}(\gamma)|^2} & \leq
& \epsilon \mu_G(B'')\sum_{\gamma \in
\wh{G}}{|\wh{f'1_{B'}}(\gamma)|^2}\\ & = &
\epsilon\mu_G(B'')\mu_G(B')\|f'\|_{L^2(\beta')}^2,
\end{eqnarray*}
by the triangle inequality and Parseval's theorem. Now $
\|f'\|_{L^2(\beta')}^2 = \alpha' - \alpha'^2$ so it follows that
this last expression is at most
$\alpha''\alpha'^2\mu_G(B'')\mu_G(B')/2^4$, and hence by the
triangle inequality and (\ref{coret2}) we have
\begin{equation*}
\sum_{\gamma \in
\Lambda}{|\wh{f''1_{B''}}(2\gamma)||\wh{f'1_{B'}}(\gamma)|^2} \geq
\alpha''\alpha'^2\mu_G(B'')\mu_G(B')/2^4.
\end{equation*}
Since
\begin{equation*}
\|f''1_{B''}\|_{L^1(\mu_G)} = 2(\alpha'' - \alpha''^2)\mu_G(B'')
\end{equation*}
we have $|\wh{f''1_{B''}}(2\gamma)| \leq 2\alpha''\mu_G(B'')$ and so
\begin{equation*}
\sum_{\gamma \in \Lambda}{|\wh{f'1_{B'}}(\gamma)|^2} \geq
\alpha'^2\mu_G(B')/2^5.
\end{equation*}
We apply Proposition \ref{localbessel} to get a system
$\wt{\mathcal{B}}'''$ with
\begin{equation*}
\dim \wt{\mathcal{B}}''' \leq 2^{10}\alpha''^{-1}\alpha'^{-2} \leq
2^{13}\alpha^{-3},
\end{equation*}
such that $\wt{\mathcal{B}}'''\cap\mathcal{B}''$ has
\begin{equation*}
\mu_G(\wt{\mathcal{B}}'''\cap\mathcal{B}'') \geq
4^{-d-2^{12}\alpha^{-3}}\mu_G(\mathcal{B}'')
\end{equation*}
and
\begin{equation*}
\Lambda \subset \{\gamma:|1-(2\gamma)(x)| \leq
2^{18}\alpha^{-3}d\eta
 \textrm{ for all
} x \in \wt{B}'''_\eta\cap B''_\eta\}.
\end{equation*}
Apply Proposition \ref{ubreg} to pick $\lambda''''$ so that
\begin{equation*}
\mathcal{B}'''':=\lambda''''((\{2x:x \in \wt{B}'''_{\rho}\})_\rho
\cap(\{2x:x \in B''_{\rho}\})_\rho)
\end{equation*}
is regular and
\begin{equation*}
\frac{\alpha^3}{2^{20}(1+d)} \leq \lambda'''' <
\frac{\alpha^3}{2^{19}(1+d)}.
\end{equation*}
Put $\wt{\mathcal{B}}'''':=\lambda''''(\{2x:x \in
\wt{B}'''_{\rho}\})_\rho$ and
$\mathcal{B}''':=\lambda''''\mathcal{B}''$. Now
\begin{eqnarray*}
\Lambda & \subset & \{\gamma : |1-(2\gamma)(x)| \leq 1/2 \textrm{
for all } x \in \wt{B}'''_{\lambda''''}\cap B''_{\lambda''''}\}\\ &
= & \{\gamma : |1-\gamma(x)| \leq 1/2 \textrm{ for all } x \in
B''''\}.
\end{eqnarray*}
$\mathcal{B}''''$ is a subsystem of
$(\alpha'/2^{14}(1+d))\mathcal{B}'$ so we may apply Lemma
\ref{l2increment} with $c=2^{-4}$ to see that
\begin{equation*}
\|1_A \ast \beta''''\|_{L^\infty(\mu_G)} \geq \alpha'(1+2^{-7}) \geq
\alpha(1+2^{-8}).
\end{equation*}
It remains only to verify the bound on the density of
$\mathcal{B}''''$. Note that
\begin{equation*}
\|1_A \ast \beta'''_{\lambda''''} \ast \beta' - 1_A \ast
\beta'\|_{L^\infty(\mu_G)} \leq 2^4d\lambda''''\lambda'' \leq
\alpha'/2
\end{equation*}
by Lemma \ref{haar}. Whence
\begin{equation*}
1_A \ast \beta'''_{\lambda''''}\ast \beta'(x') \geq 1_A \ast
\beta'(x') - \alpha'/2 \geq \alpha/2^2.
\end{equation*}
By averaging it follows that there is some $x'' \in G$ such that $
1_A \ast \beta'''_{\lambda''''}(x'') \geq \alpha/2^2$. Since $A - A$
contains no elements of order 2 we have that $x \mapsto 2x$ is
injective when restricted to $A$; we conclude that
\begin{eqnarray*}
\mu_G(2B'''_{\lambda''''}) & =& \mu_G(2(x''+B'''_{\lambda''''}))\\
& \geq & \mu_G(2(A \cap(x''+B'''_{\lambda''''})))\\& = &\mu_G(A
\cap(x''+B'''_{\lambda''''}))\\ & \geq & \frac{\alpha}{2^2}
\mu_G(B'''_{\lambda''''})\\ & \geq & \frac{\alpha}{2^2}
\left(\frac{\lambda''''}{2}\right)^{d+2^{13}\alpha^{-3}}
4^{-d-2^{12}\alpha^{-3}} \left(\frac{\lambda'\lambda''}{2}\right)^d
\mu_G(\mathcal{B}),
\end{eqnarray*}
by Lemma \ref{bourgainsize}. The claim follows.
\end{proof}
The lemma is proved.
\end{proof}

\begin{proof}[Proof of Theorem \ref{bourgrothvar}]
We construct two sequences of Bourgain systems $\wt{\mathcal{B}}_k$
and $\mathcal{B}_k'$; we write $\wt{d}_k$ for the dimension of
$\wt{\mathcal{B}}_k$, $\mathcal{B}_{k+1}$ for the intersected system
$\wt{\mathcal{B}}_k \cap \mathcal{B}_k'$, $d_k$ for the dimension of
$\mathcal{B}_k$, $\delta_k$ for the density of $\mathcal{B}_k$,
$\beta_k$ for the measure on $\mathcal{B}_k$ and
$\alpha_k:=\|1_{A_k} \ast \beta_k\|_{L^\infty(\mu_G)}$.

For $k \leq 2^{14}\log \alpha^{-1}$ we shall show inductively that
these sequences satisfy
\begin{enumerate}
\item $\wt{d}_k \leq 2^{13}\alpha^{-3}$;
\item $\mathcal{B}_k'$ is either a dilate of $\mathcal{B}_{k-1}$ or
of $(\{2x:x \in (B_{k-1})_\rho\})_\rho$;
\item $\mathcal{B}_k$ is a regular Bourgain system;
\item $d_k \leq 2d + 2^{14}\alpha^{-3}k$;
\item $\delta_k \geq \left(\frac{\alpha}{2(1+d)}\right)^{(2^8d +
2^{36}\alpha^{-3}\log \alpha^{-1})k}\mu_G(\mathcal{B})$; \item and
$\alpha_k \geq (1+2^{-12})^k\alpha$.
\end{enumerate}
We initialize the setup with $\mathcal{B}_0=\mathcal{B}$ (or, if
preferred, $\wt{\mathcal{B}}_{-1}$ as the trivial system and
$\mathcal{B}_{-1}'=\mathcal{B}$) so that the properties are
trivially satisfied. At stage $k \leq 2^{13}\log \alpha^{-1}$ apply
Lemma \ref{itlem} to $\mathcal{B}_k$. It follows that either
\begin{equation}\label{abort}
\Lambda(1_A,1_A,1_A) \geq
\frac{\alpha_k^3}{2^5}\left(\frac{\alpha_k^3}{2^{44}(1+d_k)^3}\right)^{d_k}
\mu_{G}(\mathcal{B}_k)^2;
\end{equation}
or there is a (possibly trivial) Bourgain system
$\wt{\mathcal{B}}_k$ with dimension $\wt{d}_k \leq
2^{13}\alpha_k^{-3}\leq 2^{13}\alpha^{-3}$ and another
$\mathcal{B}_k'$ which is either a dilate of $\mathcal{B}_k$ or of
$(\{2x:x \in (B_k)_\rho\})_\rho$ such that
$\mathcal{B}_{k+1}=\wt{\mathcal{B}}_k\cap\mathcal{B}_k'$ is regular
\begin{eqnarray*}
\delta_{k+1} & \geq &
\frac{\alpha_k}{2^2}\left(\frac{\alpha_k^3}{2^{22}(1+d_k)}\right)^{2^{13}\alpha_k^{-3}}
\left(\frac{\alpha_k^5}{2^{48}(1+d_k)^3}\right)^{d_k}\delta_k\\ &
\geq & \left(\frac{\alpha}{2(1+d)}\right)^{(2^8d +
2^{36}\alpha^{-3}\log \alpha^{-1})(k+1)}\mu_G(\mathcal{B})
\end{eqnarray*}
and
\begin{equation*}
\alpha_{k+1} \geq (1+2^{-12})\alpha_k \geq (1+2^{-12})^{k+1}\alpha.
\end{equation*}
It remains to check the bound on $d_{k+1}$, which follows by Lemma
\ref{join} on noting that $\mathcal{B}_{k+1}$ is the intersection of
a system of dimension $d$ and $k+1$ systems of dimension at most
$2^{13}\alpha^{-3}$.

In view of the lower bound on $\alpha_k$ and the fact that $\alpha_k
\leq 1$ it follows that there is some $k \leq 2^{13}\log
\alpha^{-1}$ such that (\ref{abort}) happens; this yields the
result.
\end{proof}

\section{An argument of Bogolio{\`u}boff and Chang}\label{bogchang}

In this section we shall prove the following proposition which draws
on techniques of Bogolio{\`u}boff \cite{NNB} as refined by Chang
\cite{MCC}. An argument of this type is contained in \cite{BJGTS2}.
\begin{proposition}\label{bg} Suppose that $G$ is a finite abelian group. Suppose that $A \subset G$ has
density $\alpha>0$ and that $|A+A| \leq K|A|$. Then there is a
regular Bourgain system $\mathcal{B}$ with
\begin{equation*}
\dim \mathcal{B} \leq 2^5K \log \alpha^{-1} \textrm{ and }
\mu_G(\mathcal{B}) \geq \left(\frac{1}{2^{14}K^2(1+\log
\alpha^{-1})}\right)^{2^4K\log \alpha^{-1}}
\end{equation*}
such that
\begin{equation*}
\|1_A \ast \beta\|_{L^\infty(\mu_G)} \geq 1/2K.
\end{equation*}
\end{proposition}
We require Chang's theorem:
\begin{proposition}\label{chng}\emph{(Chang's theorem, \cite[Proposition
3.2]{BJGIZR}).} Suppose that $A \subset G$ is a set of density
$\alpha>0$ and $\epsilon \in (0,1]$ is a parameter. Let
$\Lambda:=\{\gamma \in \wh{G}: |\wh{1_A}(\gamma)| \geq \epsilon
\alpha\}$. Then there is a set of characters $\Gamma$ with $|\Gamma|
\leq 2\epsilon^{-2}\log \alpha^{-1}$ such that $\Lambda \subset
\langle \Gamma \rangle$, where we recall that $\langle \Gamma
\rangle :=\{\sum_{\lambda \in \Gamma}{\sigma_\lambda\lambda}:\sigma
\in \{-1,0,1\}^\Gamma\}$.
\end{proposition}
\begin{proof}[Proof of Proposition \ref{bg}]
Let $\epsilon$ be a parameter to be chosen later. Apply Chang's
theorem (Proposition \ref{chng}) to the set $A$ with parameter
$\sqrt{\epsilon/3}$ to get a set of characters $\Gamma$ with
$|\Gamma| \leq 6\epsilon^{-1}\log \alpha^{-1}$ and
$\Lambda:=\{\gamma:|\wh{1_A}(\gamma)| \geq \sqrt{\epsilon/3}\alpha\}
\subset \langle \Gamma \rangle$.

Write $\mathcal{B}'$ for the Bourgain system induced by
$B(\Gamma,\epsilon/2^6(1+|\Gamma|))$ and apply Proposition
\ref{ubreg} to pick $\eta \in [1/2,1)$ so that
$\mathcal{B}:=\eta\mathcal{B}'$ is regular. It follows that
$\mathcal{B}$ has dimension at most $2|\Gamma|$ and density at least
\begin{equation*}
\left(\frac{1}{4}\right)^{2|\Gamma|}\times
\left(\frac{\epsilon}{2^6(1+|\Gamma|)}\right)^{|\Gamma|} \geq
\left(\frac{\epsilon^2}{2^{12}(1+\log
\alpha^{-1})}\right)^{|\Gamma|}.
\end{equation*}

If $\lambda \in \Lambda$ then $\lambda = \sum_{\gamma \in
\Gamma}{\sigma_\gamma\gamma}$ so
\begin{eqnarray*}
|1-\lambda(h)| \leq \sum_{\gamma \in \Gamma}{|1-\gamma(h)|}
 &= &\sum_{\gamma \in
 \Gamma}{\sqrt{2(1-\cos(4\pi\|\gamma(h)\|))}}\\ &
\leq & \sum_{\gamma \in \Gamma}{4\pi\|\gamma(h)\|} \leq 4\pi
|\Gamma| \sup_{\gamma \in \Gamma}{\|\gamma(h)\|}.
\end{eqnarray*}
So if $\lambda \in \Lambda$ then
\begin{equation*}
|1-\wh{\beta}(\lambda)|\leq \sup_{h \in B}{|1-\lambda(h)|} \leq
\epsilon/3.
\end{equation*}
Hence $|\langle 1_A \ast 1_A,1_A \ast 1_A\rangle -\langle1_A \ast
1_A,1_A \ast 1_A \ast \beta \rangle|$ is at most
\begin{eqnarray*}
|\sum_{\gamma \in
\wh{G}}{|\wh{1_A}(\gamma)|^4\left(1-\overline{\wh{\beta}(\gamma)}\right)}|
& \leq & \sup_{\gamma \in
\Lambda}{|1-\wh{\beta}(\gamma)|}\sum_{\gamma
\in\wh{G}}{|\wh{1_A}(\gamma)|^4}\\ & & +2\sup_{\gamma \not \in
\Lambda}{|\wh{1_A}(\gamma)|^2} \sum_{\gamma
\in\wh{G}}{|\wh{1_A}(\gamma)|^2}\\ & \leq &
(\epsilon/3)\alpha^2\sum_{\gamma \in\wh{G}}{|\wh{1_A}(\gamma)|^2}\\
& &+2(\epsilon/3)\alpha^2\sum_{\gamma
\in\wh{G}}{|\wh{1_A}(\gamma)|^2} \leq \epsilon\alpha^3.
\end{eqnarray*}
Moreover
\begin{equation*}
\langle 1_A \ast 1_A,1_A \ast 1_A \rangle  \geq \mu_G(\supp 1_A
\ast 1_A)^{-1}\left(\int{1_A \ast 1_A d\mu_G}\right)^2 \geq
\alpha^3/K,
\end{equation*}
by the Cauchy-Schwarz inequality and the fact that $|A+A| \leq
K|A|$. It follows from the triangle inequality that if we take
$\epsilon = 1/2K$ then
\begin{eqnarray*}
\alpha^3/2K  &\leq & |\langle 1_A \ast 1_A,1_A \ast 1_A \ast
\beta_{\Gamma,\delta} \rangle|\\ & =& |\langle 1_A \ast 1_A\ast
1_{-A} ,1_A \ast \beta_{\Gamma,\delta} \rangle|\\ & \leq & \|1_A
\ast \beta_{\Gamma,\delta}\|_{L^\infty(\mu_G)}\alpha^3.
\end{eqnarray*}
Dividing by $\alpha^3$, the result is proved.
\end{proof}

\section{The main arguments}\label{mainarg}

In this section we prove the following theorem which is the real
heart of the paper.

\begin{theorem}\label{mainthm}
Suppose that $G$ is an abelian group and $A \subset G$ is finite
with $|A+A| \leq K|A|$. If $A-A$ contains no elements of order 2
then $A$ contains at least $\exp(-CK^3\log^3(1+K))|A|^2$ three-term
arithmetic progressions for some absolute positive constant $C$.
\end{theorem}

Recall that if $G$ and $G'$ are two abelian groups with subsets
$A$ and $A'$ respectively then $\phi:A \rightarrow A'$ is a
\emph{Fre{\u\i}man homomorphism} if
\begin{equation*}
a_1+a_2=a_3+a_4 \Rightarrow \phi(a_1)+\phi(a_2) = \phi(a_3) +
\phi(a_4).
\end{equation*}
If $\phi$ has an inverse which is also a homomorphism then we say
that $\phi$ is a \emph{Fre{\u\i}man isomorphism}. For us the key
property of Fre{\u\i}man isomorphisms is that if $A$ and $A'$ are
Fre{\u\i}man isomorphic then the three-term arithmetic progressions
in $A$ and $A'$ are in one-to-one correspondence. It follows that
each set has the same number of these.

To leverage the work of \S\ref{rbvar} we need $A$ to be a large
proportion of $G$. This cannot be guaranteed but the following
proposition will allow us to move $A$ to a setting where this is
true.
\begin{proposition}\label{modl}\emph{(\cite[Proposition
1.2]{BJGIZR})}. Suppose that $G$ is an abelian group and $A \subset
G$ is finite with $|A+A| \leq K|A|$. Then there is an abelian group
$G'$ with $|G'| \leq (20K)^{10K^2}|A|$ such that $A$ is Fre{\u\i}man
isomorphic to a subset of $G'$.
\end{proposition}

\begin{proof}[Proof of Theorem \ref{mainthm}] We apply Proposition
\ref{modl} to get a subset $A'$ of a group $G'$ with density at
least $(20K)^{-10K^2}$ such that $A'$ Fre{\u\i}man isomorphic to
$A$. Since $A'$ is Fre{\u\i}man isomorphic we have $|A|=|A'|$,
$|A'+A'| \leq K|A'|$ and $A'-A'$ contains no elements of order 2. We
apply Proposition \ref{bg} to get a regular Bourgain system
$\mathcal{B}$ with
\begin{equation*}
\dim \mathcal{B} \leq 2^9K^3\log (1+K) \textrm{ and }
\mu_{G'}(\mathcal{B}) \geq (2K)^{-2^{13}K^3\log (1+K)}
\end{equation*}
such that $\|1_{A'} \ast \beta\|_{L^\infty(\mu_{G'})} \geq 1/2K$. We
now apply Theorem \ref{bourgrothvar} to get the result.
\end{proof}

The proof of Theorem \ref{premainthm} is now rather straightforward.

\begin{proof}[Proof of Theorem \ref{premainthm}]
Write $K:=|A+A|/|A|$ and suppose that $a, a' \in A$ have $a-a'$ of
order 2. Then $a+a=2a'$ is a non-trivial three-term progression in
$A$ which contradicts the hypothesis. It follows that we may apply
Theorem \ref{mainthm} to conclude that $A$ contains at least
$\exp(-CK^3\log^3(1+K))|A|^2$ progressions; however we know this to
be at most $|A|$, whence
\begin{equation*}
\exp(CK^3\log^3(1+K)) \geq |A|.
\end{equation*}
The result follows on rearranging.
\end{proof}

Proving Theorem \ref{schoent} simply requires us to apply Theorem
\ref{premainthm} in moreorless the same manner as Schoen applies
Theorem \ref{RuzBou}.
\begin{proof}[Proof of Theorem \ref{schoent}] Write
\begin{equation*}
  S:=\{a \in A \cap B: \not \exists a' \in A, b' \in B \textrm{ with } a' \neq b' \textrm{ such
  that } a'+b' =2a\},
\end{equation*}
and note that crucially we have
\begin{equation}\label{ky}
(A+B) \setminus (A\ \wh{+}\ B) = 2S,
\end{equation}
and moreover that $S$ contains no three-term progressions $(a,b,c)$
with $a+b=2c$ and $a \neq b$.

Let $S'$ be a subset of $S$ such that for all $s \in 2S$ there is
exactly one $s' \in S'$ such that $2s' = s$. It is easy to see that
$|S'| =|2S|$.

We claim that $S'$ contains no non-trivial three-term progressions.
Suppose that $a,b,c \in S'$ have $a+b=2c$. Since $S' \subset S$ we
conclude that $a=b$, but in this case we have $2a=2c$ which, by
choice of $S'$, implies that $a=c$. The claim follows.

Consequently we may apply Theorem \ref{premainthm} to conclude that
\begin{equation*}
|S'+S'| \gg |S'| \left(\frac{\log |S'|}{(\log \log
|S'|)^3}\right)^{\frac{1}{3}}.
\end{equation*}
Recalling that $n=|A+B|$ we can rearrange this expression to give
\begin{equation*}
|S'| \ll |S'+S'|\left(\frac{(\log \log |S'+S'|)^3}{\log
|S'+S'|}\right)^{\frac{1}{3}} \ll |A+B|\left(\frac{(\log \log
n)^3}{\log n}\right)^{\frac{1}{3}},
\end{equation*}
since the middle expression is an increasing function of $|S'+S'|$
and $S'+S' \subset A+B$. The result follows from (\ref{ky}) and the
fact that $|S'|=|2S|$.
\end{proof}

\section{Concluding remarks}\label{concluding remarks}

The extension of Theorem \ref{RuzBou} to the groups $\Z^r$ and
$\Z/N\Z$ (with the same bound) is implicit in the works of Ruzsa
\cite{IZR3AP}, Bourgain \cite{JB} and Stanchescu \cite{YVS}.
Moreover, since there are particularly good versions of the modelling
proposition (Proposition \ref{modl}) for these groups it seems very
likely that our Proposition \ref{bg} could be used in conjunction
with a more traditional $\ell^\infty$-density increment argument of
Bourgain \cite{JB} to prove the following.
\begin{theorem}
Suppose that $G$ is $\Z^r$ or $\Z/N\Z$ and $A \subset G$ is finite
with $|A+A| \leq K|A|$. Then $A$ contains at least
$\exp(-CK^{2+o(1)} )|A|^2$ three-term arithmetic progressions for
some absolute $C>0$.
\end{theorem}
Indeed, it appears that with the methods of \cite{TSFT} one could
replace $K^{2+o(1)}$ by $K^2 \log (1+K)$, thereby directly
generalizing Bourgain's version of Roth's theorem from \cite{JB}.

We have not considered how the ideas in Bourgain's recent paper
\cite{JBRoth2} might come into play to give an even stronger result;
the following is a natural question.
\begin{problem}
Find a direct generalization of the result of \cite{JBRoth2} to sets
with small sumset. That is, show that if $A \subset \Z/N\Z$ is
finite with $|A+A| \leq K|A|$, then $A$ contains at least
$\exp(-CK^{3/2}\log^3(1+K) )|A|^2$ three-term arithmetic
progressions for some absolute $C>0$.
\end{problem}
Among other things Theorem \ref{premainthm} immediately improves a
result of Stanchescu \cite{YVS} who, answering a further question of
Fre{\u\i}man \cite{GAF}, used Theorem \ref{RuzBou} to bound from
below the size of $|A+A|/|A|$ when $A \subset \Z^2$ is finite and
contains no three collinear points. This is an intriguing question
because one appears to have so much extra information to play with:
not only does $A$ not contain any three-term progressions but it
also avoids any triples $(a,b,c)$ with $\lambda a+\mu b=(\lambda +
\mu)c$ for any positive integers $\lambda$ and $\mu$.
\begin{problem}
Find a constant an absolute constant $c>2/3$ such that if $A \subset
\Z^2$ is finite and contains no three collinear points then $|A+A|
\gg |A|\log^c |A|$.
\end{problem}

Moves to generalize additive problems to arbitrary abelian groups
have also spawned the observation (see, for example, \cite{BJGFFM}
and \cite{BJGTCTI}) that some arguments can be modelled very cleanly
(and often more effectively) in certain well behaved abelian groups.
Again, it would be surprising if one could not prove the following
using the methods outlined above.
\begin{theorem}
Suppose that $G$ is a vector space over $\mathbb{F}_3$ and $A
\subset G$ is finite with $|A+A| \leq K|A|$. Then $A$ contains at
least $\exp(-CK)|A|^2$ three-term arithmetic progressions for some
absolute constant $C>0$.
\end{theorem}

In a different direction it maybe that the following problem
captures the essence of Roth's theorem in a natural general setting.
\begin{problem}
Suppose that $A \subset \Z$ has at least $\delta|A|^3$ additive
quadruples. Find a good absolute constant $c>0$ such that we can
conclude that $A$ contains at least $\exp(-C\delta^{-c})|A|^2$ three
term arithmetic progressions.
\end{problem}
It is immediate from the quantitative Balog-Szemer{\'e}di-Gowers
theorem (see the paper \cite{WTG} of Gowers) that there is some
$c>0$; the problem is to find a good value.

\section{Acknowledgments}
I should like to thank Tim Gowers and Ben Green for supervision, and
Ben Green and Terry Tao for making the preprint \cite{BJGTCTF}
available.

\bibliographystyle{alpha}

\bibliography{master}

\end{document}